\documentclass[11pt]{article}
\usepackage{amssymb}
\usepackage{amsmath}
\usepackage{fancyhdr}
\usepackage{amsthm}
\usepackage[T1]{fontenc}
\usepackage{afterpage}  %

\usepackage {amsfonts, amsmath, array}

\usepackage [mathscr] {eucal}
\usepackage {amscd}
\usepackage [T2A] {fontenc}
\usepackage {graphics}
\usepackage [dvips] {color}
\usepackage [cp866] {inputenc}

\theoremstyle{plain}

\newtheorem{example}{Example}

\title{On prolongations of quasigroups}
\author{Ivan I. Deriyenko and Wieslaw A.
Dudek}
\date{}

\begin{document}
\sloppy \maketitle

\begin{abstract}
We prove that any quasigroup admissing complete or quasicomplete
mapping has a prolongation to a quasigroup having one element
more.
\end{abstract}

\section{Introduction}

By a {\it prolongation} of a quasigroup we mean a process which
shows how, starting from a quasigroup $Q(\cdot)$ of order $n$, we
can obtain a quasigroup $Q^{\prime}(\circ)$ of order $n+1$ such
that the set $Q^{\prime}$ is obtained from the set $Q$ by the
adjunction of one additional element. In other words, it is a
process which shows how a given Latin square extends to a new
Latin square by the adjunction of one additional row and one
column. The first construction of prolongation was proposed by R.
H. Bruck \cite{Br} who considered only the case of idempotent
quasigroups. More general construction was given by J. D\'enes and
K. P\'asztor \cite{DP}. Further generalizations, for special types
of quasigroups, have been discussed in \cite{VD1} and \cite{VD2}
by V. D. Belousov. In fact, the construction proposed by V. D.
Belousov is more elegant form of the construction proposed by J.
D\'enes and K. P\'asztor. G. B. Belyavskaya studied this problem
together with the inverse problem, i.e., with the problem how from
a given Latin square of order $n$ one can obtain a Latin square of
order $n-1$ (cf. \cite{GB1, GB2, GB3}). Quasigroups obtained by
the construction proposed by G. B. Belyavskaya are not isotopic to
quasigroups obtained by the constructions proposed by R. H. Bruck
and V.D. Belousov. This means that we have two different methods
of construction of prolongations.

Below we present a third method. Our method can be applied to any
quasigroup of order $n$ with the property that its multiplication
table possesses a partial transversal of length $n-1$, i.e., a
sequence of $n-1$ distinct elements contained in distinct rows and
distinct columns. All these three constructions are presented in
short elegant form.

\section{Definitions and basic facts}

In this paper $Q(\cdot)$ always denotes a quasigroup. The set
$Q^{\prime}$ is identified with the set $Q\cup\{q\}$, where
$q\notin Q$.

Any mapping $\sigma$ of a quasigroup $Q(\cdot)$ defines on $Q$ a
new mapping $\overline{\sigma}$, called {\it conjugated to}
$\sigma$, such that
\begin{equation}\label{1}
\overline{\sigma}(x)=x\cdot\sigma(x)
\end{equation}
for all $x\in Q$. If $\sigma$ is the identity mapping
$\varepsilon$, then $\overline{\sigma}(x)=x^2$. The set
$${\rm def}(\sigma)=Q\setminus\overline{\sigma}(Q),$$
where $\overline{\sigma}(Q)=\{\overline{\sigma}(x)\,|\,x\in Q\}$,
is called the {\it defect} of $\sigma$.

A mapping $\sigma$ is {\it quasicomplete} on a quasigroup
$Q(\cdot)$ if $\overline{\sigma}(Q)$ contains all elements of $Q$
except one. In this case there exists an element $a\in Q$, called
{\it special}, such that
$a=\overline{\sigma}(x_{1})=\overline{\sigma}(x_{2})$ for some
$x_{1},x_{2}\in Q$. If $\overline{\sigma}(Q)$ contains all
elements of $Q$, then we say that $\sigma$ is {\it complete}. A
quasigroup having at least one complete mapping is called {\it
admissible}. V. D. Belousov proved in \cite{VD2} (see also
\cite{VD1}) that any admissible quasigroup is isotopic to some
idempotent quasigroup and has a prolongation. Since for a given
admissible quasigroup the method of constructions of a
prolongation proposed by V.D. Belousov gives, in fact, a
quasigroup which is isotopic to a quasigroup obtained from the
corresponding idempotent quasigroup (by the method proposed by R.
H. Bruck) we will identify these two methods and will call it the
{\it classical construction}.

\section{Prolongations of admissible quasigroups}

{\bf 1. Classical constructions.} The idea of the construction
proposed by R. H. Bruck is presented by the following tables,
where the corresponding empty cells of these tables are identical.
\[
\begin{array}{rcl}
\begin{array}{c|cccccc|}
\cdot&1&2&3&4&\ldots&n\\ \hline
 1&1&&&&&\\
 2&&2&&&&\\
 3&&&3&&&\\
 4&&&&4&&\\
 \vdots&&&&&\ddots&\\
 n&&&&&&n\\ \hline
\end{array}
&\;\;\;\longrightarrow\;\;\;&
\begin{array}{c|cccccc|c}
\circ&1&2&3&4&\ldots&n&q\\ \hline
 1&q&&&&&&1\\
 2&&q&&&&&2\\
 3&&&q&&&&3\\
 4&&&&q&&&4\\
 \vdots&&&&&\ddots&&\vdots\\
 n&&&&&&q&n\\ \hline
 q&1&2&3&4&\ldots&n&q\\
\end{array}
\end{array}
\]

\bigskip

The quasigroup $Q^{\prime}(\circ)$ obtained from the quasigroup
$Q(\cdot)$ is a loop with the identity $q$. The operation on
$Q^{\prime}$ is defined according to the formula:
\begin{equation}\label{e2}
x\circ y=\left\{\begin{array}{clll}
 x\cdot y&{\rm for}& x,y\in Q,\ x\ne y,\\
 x&{\rm for}&x\in Q,\ y=q,\\
 y&{\rm for}&x=q,\ y\in Q,\\
 q&{\rm for}& x=y\in Q^{\prime}.
\end{array}\right.
\end{equation}

In the construction for a prolongation of an admissible quasigroup
$Q(\cdot)$ proposed by V. D. Belousov \cite{VD2} the complete
mapping $\sigma$ of $Q(\cdot)$ and its conjugated mapping
$\overline{\sigma}$ are used. The operation on $Q^{\prime}$ is
defined by the formula:
\begin{equation}\label{e3}
 x\circ y=\left\{\begin{array}{clll}
 x\cdot y&{\rm for}&x,y\in Q,\ y\ne\sigma(x),\\[3pt]
 \overline{\sigma}(x)&{\rm for}&x\in Q,\ y=q,\\[3pt]
 \overline{\sigma}\sigma^{-1}\!(y)&{\rm for}&x=q,\ y\in Q,\\[3pt]
 q&{\rm for}&x\in Q,\ y=\sigma(x),\\[3pt]
 q&{\rm for}& x=y=q.
\end{array}\right.
\end{equation}

Geometrically this means that the multiplication table (Latin
square) $L^{\prime}=[a'_{ij}]$ of a quasigroup $Q^{\prime}(\circ)$
is obtained from the multiplication table $L=[a_{ij}]$ of a
quasigroup $Q(\cdot)$ by the adjunction of one row and one column
in this way that all elements from the cells $a_{i\sigma(i)}$ are
moved to the last place of the $i$-th row and $\sigma(i)$-th
column of $L^{\prime}$. Elements of the cells $a_{i\sigma(i)}$ are
replaced by $q=n+1$. Additionally we put $a_{qq}=q$. In other
words: $a'_{ij}=a_{ij}$ for $i\ne\sigma(i)$,\
$a'_{iq}=a_{i\sigma(i)}=\overline{\sigma}(i)$,\
$a'_{qj}=a_{\sigma^{-1}(j)j}=\overline{\sigma}\sigma^{-1}(j)$\
and\ $a'_{i\sigma(i)}=a'_{qq}=q$.

\begin{example}\label{ex1}\rm
Consider the quasigroup $Q(\cdot)$ with the multiplication table
\[\begin{array}{c|ccccc}
\cdot&1&2&3&4&5\\ \hline\rule{0mm}{4mm}
 1&1&2&3&4&5\\
 2&4&3&1&5&2\\
 3&2&5&4&1&3\\
 4&5&4&2&3&1\\
 5&3&1&5&2&4
\end{array}
\]
and its two complete mappings
$\sigma=\left(\arraycolsep=.5mm\begin{array}{ccccccc}
1&2&3&4&5 \\
4&2&1&5&3 \end{array}\right)$ and
$\tau=\left(\arraycolsep=.5mm\begin{array}{ccccccc}
1&2&3&4&5 \\
3&1&2&5&4 \end{array}\right)$. Then, as it is not difficult to
see,
$\overline{\sigma}=\left(\arraycolsep=.5mm\begin{array}{ccccccc}
1&2&3&4&5 \\
4&3&2&1&5 \end{array}\right)$ and
$\overline{\tau}=\left(\arraycolsep=.5mm\begin{array}{ccccccc}
1&2&3&4&5 \\
3&4&5&1&2 \end{array}\right)$. Using these two mappings we can
construct two different prolongations:
\[\begin{array}{rcl}\begin{array}{c|ccccc|c}
\circ_1&1&2&3&4&5&6\\ \hline\rule{0mm}{4mm}
 1&1&2&3&6&5&4\\
 2&4&6&1&5&2&3\\
 3&6&5&4&1&3&2\\
 4&5&4&2&3&6&1\\
 5&3&1&6&2&4&5\\ \hline\rule{0mm}{4mm}
 6&2&3&5&4&1&6
\end{array}
 & \ \ \ \ &
 \begin{array}{c|ccccc|c}
\circ_2&1&2&3&4&5&6\\ \hline\rule{0mm}{4mm}
 1&1&2&6&4&5&3\\
 2&6&3&1&5&2&4\\
 3&2&6&4&1&3&5\\
 4&5&4&2&3&6&1\\
 5&3&1&5&6&4&2\\ \hline\rule{0mm}{4mm}
 6&4&5&3&2&1&6\\
\end{array}
\end{array}
\]
The first prolongation is obtained by $\sigma$, the second by
$\tau$.

\medskip

By transpositions of rows and columns, we can transform these two
tables into multiplication tables of loops $Q^{\prime}(\star_1)$
and $Q^{\prime}(\star_2)$:

\[\begin{array}{rcl}\begin{array}{c|ccccc|c}
\star_1&1&2&3&4&5&6\\ \hline\rule{0mm}{4mm}
 1&1&2&3&4&5&6\\
 2&2&3&5&6&1&4\\
 3&3&1&6&5&4&2\\
 4&4&6&1&3&2&5\\
 5&5&4&2&1&6&3\\ \hline\rule{0mm}{4mm}
 6&6&5&4&2&5&1
\end{array}
 & \ \ \ \ &
 \begin{array}{c|ccccc|c}
\star_2&1&2&3&4&5&6\\ \hline\rule{0mm}{4mm}
 1&1&2&3&4&5&6\\
 2&2&6&5&1&3&4\\
 3&3&1&2&6&4&5\\
 4&4&5&6&2&1&3\\
 5&5&4&1&3&6&2\\ \hline\rule{0mm}{4mm}
 6&6&3&4&5&2&1
\end{array}
\end{array}
\]

Since $\gamma(x\star_1y)=\alpha(x)\star_2\beta(y)$, where
\[
\alpha=\left(\arraycolsep=.5mm\begin{array}{cccccccc}
1&2&3&4&5&6 \\
2&4&1&6&5&3 \end{array}\right), \ \ \ \
\beta=\left(\arraycolsep=.5mm\begin{array}{cccccccc}
1&2&3&4&5&6 \\
4&1&6&3&5&2 \end{array}\right), \ \ \ \
\gamma=\left(\arraycolsep=.5mm\begin{array}{cccccccc}
1&2&3&4&5&6 \\
1&2&4&5&3&6 \end{array}\right),
\]
loops $Q^{\prime}(\star_1)$ and $Q^{\prime}(\star_2)$ are
isotopic. This means that also prolongations $Q^{\prime}(\circ_1)$
and $Q^{\prime}(\circ_2)$ are isotopic. \hfill$\Box{}$
\end{example}

\medskip

If the diagonal of the multiplication table of a quasigroup
$Q(\cdot)$ contains all elements of $Q$, then as $\sigma$ we can
select the identity mapping. In this case the formula \eqref{e3}
has the form:
\begin{equation}\label{e4}
 x\circ y=\left\{\begin{array}{clll}
 x\cdot y&{\rm for}&x,y\in Q,\ x\ne y,\\
 x^2&{\rm for}&x\in Q,\ y=q,\\
 y^2&{\rm for}&x=q,\ y\in Q,\\
 q&{\rm for}& x=y\in Q^{\prime}.
\end{array}\right.
\end{equation}

If $(Q(\cdot)$ is an idempotent quasigroup, then \eqref{e4}
coincides with \eqref{e2} and $Q^{\prime}(\circ)$ is a loop with
the identity $q$.

\begin{example}\label{ex2}\rm
The diagonal of the multiplication table of the additive group
$\mathbb{Z}_3$ contains all elements of $\mathbb{Z}_3$. So,
according to \eqref{e4}, the prolongation
$\mathbb{Z}^{\prime}_3(\circ)$ has the following multiplication
table:
\[
\begin{array}{c|ccc|c}
\circ&0&1&2&3\\ \hline\rule{0mm}{4mm}
      0&3&1&2&0\\
      1&1&3&0&2\\
      2&2&0&3&1\\ \hline\rule{0mm}{4mm}
      3&0&2&1&3
\end{array}
\]
Putting $\alpha=\left(\arraycolsep=.5mm\begin{array}{cccc}
0&1&2&3 \\
3&1&2&0 \end{array}\right)$ and $x\odot y=\alpha(x\circ y)$ we can
see that $\mathbb{Z}^{\prime}_3(\circ)$ is isotopic to the Klein's
group $K_4(\odot)$.

Using $\sigma=\left(\arraycolsep=.5mm\begin{array}{ccc}
0&1&2 \\
2&0&1 \end{array}\right)$ and
$\tau=\left(\arraycolsep=.5mm\begin{array}{cccc}
0&1&2 \\
1&2&0 \end{array}\right)$ we obtain two non-commutative
prolongations:
\[
\begin{array}{ccc}
\begin{array}{c|ccc|c}
\circ&0&1&2&3\\ \hline\rule{0mm}{4mm}
    0&0&1&3&2\\
    1&3&2&0&1\\
    2&2&3&1&0\\ \hline\rule{0mm}{4mm}
    3&1&0&2&3
\end{array}
& \rule{10mm}{0mm} &
\begin{array}{c|ccc|c}
\circ&0&1&2&3\\ \hline\rule{0mm}{4mm}
    0&0&3&2&1\\
    1&1&2&3&0\\
    2&3&0&1&2\\ \hline\rule{0mm}{4mm}
    3&2&1&0&3
\end{array}
\end{array}
\]
These prolongations also are isotopic to the Klein's group. For
the first we have $x\odot y=\alpha(x)\circ\beta(y)$, for the
second $x\odot y=\beta(x)\circ\alpha(y)$, where
$\alpha=\left(\arraycolsep=.5mm\begin{array}{cccc}
0&1&2&3 \\
0&3&2&1 \end{array}\right)$ and
$\beta=\left(\arraycolsep=.5mm\begin{array}{cccc}
0&1&2&3 \\
0&1&3&2 \end{array}\right)$. \hfill$\Box{}$
\end{example}

\noindent {\bf 2. The construction proposed by G. B. Belyavskaya.}
This construction is valid for admissible quasigroups. At first we
consider the case when $Q(\cdot)$ is an idempotent quasigroup. To
find the prolongation $Q^{\prime}(\diamond)$ of $Q(\cdot)$ we
select an arbitrary element $a\in Q$. Next, in the multiplication
table of $Q(\cdot)$ we replace all elements of the diagonal,
except $a$, by $q$ and adjunct one column and one row:

\[
\arraycolsep=.9mm
\begin{array}{rcl}
\begin{array}{c|cccccc|}
\cdot&1&2&\ldots&a&\ldots&n\\ \hline\rule{0mm}{4mm}
 1&1&&&&&\\
 2&&2&&&&\\
  \vdots&&&\ddots&&&\\
  a&&&&a&&\\
 \vdots&&&&&\ddots&\\
 n&&&&&&n\\ \hline
\end{array}
&\;\;\;\longrightarrow\;\;\;&\arraycolsep=1mm
\begin{array}{c|cccccc|c}
\diamond&1&2&\ldots&a&\ldots&n&q\\ \hline\rule{0mm}{4mm}
 1&q&&&&&&1\\
 2&&q&&&&&2\\
 \vdots&&&\ddots&&&&\vdots\\
 a&&&&a&&&q\\
 \vdots&&&&&\ddots&&\vdots\\
 n&&&&&&q&n\\ \hline
 q&1&2&\ldots&q&\ldots&n&a\\
\end{array}
\end{array}
\]

\bigskip

The operation in $Q(\diamond)$ is defined in the following way:

\begin{equation}\label{e5}
 x\diamond y=\left\{\begin{array}{clll}
 x\cdot y&{\rm for}&x,y\in Q,\ x\ne y,\\
 q&{\rm for}&x=y\in Q-\{a\},\\
 a&{\rm for}&x=y=a,\\
 x&{\rm for}&x\in Q-\{a\},\ y=q,\\
 y&{\rm for}&x=q,\ y\in Q-\{a\},\\
 q&{\rm for}&x=q,\ y=a,\\
 q&{\rm for}&x=a,\ y=q,\\
 a&{\rm for}&x=y=q.\end{array}\right.
\end{equation}

In a general case, when $Q(\cdot)$ is ''only" an admissible
quasigroup, we can select a complete mapping $\sigma$ of $Q$ and
fix an arbitrary element $a\in Q$. Then, obviously, there exists
an uniquely determined element $x_a\in Q$ such that
$a=x_a\cdot\sigma(x_a)$. The prolongation $Q^{\prime}(\diamond)$
of $Q(\cdot)$ can be defined by

\begin{equation}\label{e6}
 x\diamond y=\left\{\begin{array}{clll}
 x\cdot y&{\rm for}&x,y\in Q,\ y\ne\sigma(x),\\
 q&{\rm for}&x\in Q-\{x_a\},\ y=\sigma(x),\\
 a&{\rm for}&x=x_a,\ y=\sigma(x_a),\\
 \overline{\sigma}(x)&{\rm for}&x\in Q-\{x_a\},\ y=q,\\
 \overline{\sigma}\sigma^{-1}(y)&{\rm for}&x=q,\ y\ne\sigma(x_a),\\
 q&{\rm for}&x=q,\ y=\sigma(x_a),\\
 q&{\rm for}&x=x_a,\ y=q,\\
 a&{\rm for}&x=y=q.\end{array}\right.
\end{equation}
Selecting different $\sigma$ and different $a$ we obtain different
prolongations.

From a formal point of view, the above construction is a
generalization on the classical construction. Indeed, putting
$\sigma(q)=q$ we extend $\sigma$ to a complete mapping of
$Q^{\prime}$. Next, putting $a=x_a=q$ in \eqref{e6} we obtain
\eqref{e3}.

If the diagonal of the multiplication table of $Q(\cdot)$ contains
all elements of $Q$, then as $\sigma$ can be selected the identity
mapping and the formula \eqref{e6} can be written in the form:
\begin{equation}\label{e7}
 x\diamond y=\left\{\begin{array}{clll}
 x\cdot y&{\rm for}&x,y\in Q,\ x\ne y,\\
 q&{\rm for}&x=y\in Q-\{x_a\},\\
 a&{\rm for}&x=y=x_a,\\
 x^2&{\rm for}&x\in Q-\{x_a\},\ y=q,\\
 y^2&{\rm for}&x=q,\ y\in Q-\{x_a\},\\
 q&{\rm for}&x=q,\ y=x_a,\\
 q&{\rm for}&x=x_a,\ y=q,\\
 a&{\rm for}&x=y=q.\end{array}\right.
\end{equation}

For idempotent quasigroups it coincides with \eqref{e5} but,
generally, prolongations obtained by the method proposed by G. B.
Belyavskaya are not isotopic to prolongations obtained by the
method proposed by V. D. Belousov. Below we present the
corresponding example.

\begin{example}\label{ex4}\rm The prolongation $\mathbb{Z}^{\prime}_3(\diamond)$
of the additive group $\mathbb{Z}_3$
constructed according to \eqref{e7}, where $a=1$, $x_a=2$, $q=3$,
has the following multiplication table:
\[
\begin{array}{c|ccc|c}
\diamond&0&1&2&3\\ \hline\rule{0mm}{4mm}
       0&3&1&2&0\\
       1&1&3&0&2\\
       2&2&0&1&3\\ \hline\rule{0mm}{4mm}
       3&0&2&3&1
\end{array}
\]
This prolongation is isotopic to the group $\mathbb{Z}_4(+)$. The
connection between $\mathbb{Z}_4(+)$ and
$\mathbb{Z}^{\prime}_3(\diamond)$ is given by the formula
$\gamma(x+y)=\alpha(x)\diamond\alpha(y)$, where
$\alpha=\left(\arraycolsep=.5mm\begin{array}{cccc}
0&1&2&3 \\
3&1&2&0 \end{array}\right)$, \
$\gamma=\left(\arraycolsep=.5mm\begin{array}{cccc}
0&1&2&3 \\
1&2&3&0 \end{array}\right)$. So, the prolongation of
$\mathbb{Z}_3$ constructed by \eqref{e7} and the prolongation of
$\mathbb{Z}_3$ constructed by \eqref{e4} (in Example \ref{ex2})
are not isotopic. \hfill$\Box{}$
\end{example}

\begin{example}\label{ex3}\rm Let $Q(\cdot)$ and $\sigma$ be as in Example \ref{ex1}.
Then, for example, for $a=2$ we have $x_a=3$. Whence, according to
\eqref{e6}, we obtain the prolongation:

\[
\begin{array}{c|ccccc|c}
\diamond&1&2&3&4&5&6\\ \hline\rule{0mm}{4mm}
       1&1&2&3&6&5&4\\
       2&4&6&1&5&2&3\\
       3&2&5&4&1&3&6\\
       4&5&4&2&3&6&1\\
       5&3&1&6&2&4&5\\ \hline\rule{0mm}{4mm}
       6&6&3&5&4&1&2
\end{array}
\]

Similarly, for $a=3$ we have $x_a=2$ and consequently
\[
\begin{array}{c|ccccc|c}
\diamond&1&2&3&4&5&6\\ \hline\rule{0mm}{4mm}
       1&1&2&3&6&5&4\\
       2&4&3&1&5&2&6\\
       3&6&5&4&1&3&2\\
       4&5&4&2&3&6&1\\
       5&3&1&6&2&4&5\\ \hline\rule{0mm}{4mm}
       6&2&6&5&4&1&3
\end{array}
\]
Applying Theorem 2.5 from \cite{der3} we can verify that these
prolongations are not isotopic to the prolongation obtained in
Example \ref{ex1}. \hfill$\Box{}$
\end{example}

\section{Our construction}

In the previous section methods of construction of a prolongation
of quasigroups that have a complete mapping were given. But, as it
is proved in \cite{Hall} (see also \cite{DK}, p. 36) there are
quasigroups which do not possess such mappings. For example, a
group of order $4k+2$ has no complete mapping.

Below, we give a new method of a construction of prolongations for
quasigroups that have a quasicomplete mapping. Our method can also
be applied to quasigroups that have a complete mapping.

Let $Q(\cdot)$ be an arbitrary quasigroup with a quasicomplete
mapping $\sigma$. Then $|\overline{\sigma}(Q)|=n-1$ and
def$(\sigma)=d$ for some $d\in Q$. In this case we also have
$\overline{\sigma}(x_1)=\overline{\sigma}(x_2)=a$, i.e.,
$x_1\cdot\sigma(x_1)=x_2\cdot\sigma(x_2)=a$ in $Q(\cdot)$, for
some $x_1,x_2,a\in Q$, $x_1\ne x_2$.

The idea of our construction is presented by the following tables,
where for simplicity it is assumed that $\sigma$ is the identity
mapping and all elements of $Q$, except $x_1$ and $x_2$, are
idempotents.

\[
\arraycolsep=.9mm
\begin{array}{rcl}
\begin{array}{c|cccccccc|}
\cdot&1&2&\ldots&x_1&\ldots&x_2&\ldots&n\\ \hline
 1&1&&&&&&&\\
 2&&2&&&&&&\\
  \vdots&&&\ddots&&&&&\\
 x_1&&&&a&&&&\\
 \vdots&&&&&\ddots&&&\\
 x_2&&&&&&a&&\\
 \vdots&&&&&&&\ddots&\\
 n&&&&&&&&n\\ \hline
\end{array}
&\;\;\;\longrightarrow\;\;\;&\arraycolsep=1mm
\begin{array}{c|cccccccc|c}
\ast&1&2&\ldots&x_1&\ldots&x_2&\ldots&n&q\\ \hline
 1&q&&&&&&&&1\\
 2&&q&&&&&&&2\\
  \vdots&&&\ddots&&&&&&\vdots\\
  x_1&&&&a&&&&&q\\
  \vdots&&&&&\ddots&&&&\vdots\\
  x_2&&&&&&q&&&a\\
 \vdots&&&&&&&\ddots&&\vdots\\
 n&&&&&&&&q&n\\ \hline
 q&1&2&\ldots&q&\ldots&a&\ldots&n&d\\
\end{array}
\end{array}
\]

This new table is obtained from the old one by replacing all
elements of the diagonal, except $a=x_1\cdot x_1$, by $q$ and
adding one new row and column such that $x\ast q=q\ast x=x$ for
$x\in Q-\{x_1,x_2\}$, $x_1\ast q=q\ast x_1=q$, $x_2\ast q=q\ast
x_2=a$, $q\ast q=d$.

The operation of this new quasigroup is determined by the formula:
\begin{equation}\label{e8}
 x\ast y=\left\{\begin{array}{clll}
 x\cdot y&{\rm for}&x,y\in Q,\ x\ne y,\\
  q&{\rm for}&x=y\in Q-\{x_1\},\\
 a&{\rm for}&x=y=x_1,\\
 x&{\rm for}&x\in Q-\{x_1,x_2\},\ y=q,\\
 y&{\rm for}&x=q,\ y\in Q-\{x_1,x_2\},\\
 q&{\rm for}&x=x_1,\ y=q \ \ or \ \ x=q,\ y=x_1,\\
 a&{\rm for}&x=x_2,\ y=q \ \ or \ \ x=q,\ y=x_2,\\
 d&{\rm for}&x=y=q.\end{array}\right.
\end{equation}

In the general case, when $\sigma$ is an arbitrary quasicomplete
mapping of $Q$, def$(\sigma)=d$,
$a=\overline{\sigma}(x_1)=\overline{\sigma}(x_2)$, $x_1\ne x_2$
and $x_1$ is fixed, the operation of $Q^{\prime}(\ast)$ has the
form:
\begin{equation}\label{e9}
 x\ast y=\left\{\begin{array}{clll}
 x\cdot y&{\rm for}&x,y\in Q,\ y\ne\sigma(x),\\
  q&{\rm for}&x\in Q-\{x_1\},\ y=\sigma(x),\\
 a&{\rm for}&x=x_1,\ y=\sigma(x),\\
 \overline{\sigma}(x)&{\rm for}&x\in Q-\{x_1,x_2\},\ y=q,\\
 \overline{\sigma}\sigma^{-1}(y)&{\rm for}&x=q,\ y\ne\sigma(x_1),\ y\ne\sigma(x_2),\\
 q&{\rm for}&x=x_1,\ y=q \ \ or \ \ x=q,\ y=\sigma(x_1),\\
 a&{\rm for}&x=x_2,\ y=q \ \ or \ \ x=q,\ y=\sigma(x_2),\\
 d&{\rm for}&x=y=q.\end{array}\right.
\end{equation}

If in the above formula we delete $x_2$ and assume that $\sigma$
is a complete mapping, then for $x_1=x_a$ and $d=a$ this formula
will be identical with \eqref{e7}. This means that our
construction is a generalization of the construction proposed by
G. B. Belyavskaya. Consequently, it is also a generalization of
the classical construction.

\begin{example}\label{ex5}\rm
Let $Q(\cdot)$ be a quasigroup defined in Example \ref{ex1}. The
mapping $\sigma=\left(\arraycolsep=.5mm\begin{array}{ccccccc}
1&2&3&4&5 \\
4&5&2&3&1 \end{array}\right)$ is quasicomplete on $Q$,
$\overline{\sigma}=\left(\arraycolsep=.5mm\begin{array}{ccccccc}
1&2&3&4&5 \\
4&2&5&2&3 \end{array}\right)$ is its conjugated mapping,
def$(\sigma)=1$, $\overline{\sigma}(2)=\overline{\sigma}(4)=2$.
Hence $d=1$, $a=2$, $x_1=2$, $x_2=4$. Putting $q=6$ and using our
construction we obtain the following prolongation of $Q(\cdot)$:
\[\begin{array}{c|ccccc|c}
 \ast&1&2&3&4&5&6\\ \hline\rule{0mm}{4mm}
    1&1&2&3&6&5&4\\
    2&4&3&1&5&2&6\\
    3&2&6&4&1&3&5\\
    4&5&4&6&3&1&2\\
    5&6&1&5&2&4&3\\ \hline\rule{0mm}{4mm}
    6&3&5&2&4&6&1
    \end{array}
\]

For $x_1=4$, $x_2=2$ our construction gives the quasigroup:
\[\begin{array}{c|ccccc|c}
 \ast&1&2&3&4&5&6\\ \hline\rule{0mm}{4mm}
    1&1&2&3&6&5&4\\
    2&4&3&1&5&6&2\\
    3&2&6&4&1&3&5\\
    4&5&4&2&3&1&6\\
    5&6&1&5&2&4&3\\ \hline\rule{0mm}{4mm}
    6&3&5&6&4&2&1
    \end{array}
\]
From Theorem 2.5 in \cite{der3} it follows that these two
prolongations are isotopic, but they are not isotopic to the
prolongation constructed in Example~\ref{ex1} and in Example
\ref{ex3}. \hfill$\Box{}$
\end{example}

\section{Conclusion}

The Brualdi conjecture (cf. \cite{DK}, p.103) says that each Latin
square $n\times n$ possesses a sequence of $k\geqslant n-1$
distinct elements selected from different rows and different
columns. In other words, each finite quasigroup has at least one
complete or quasicomplete mapping. It is known that if a
quasigroup $Q(\cdot)$ has a complete mapping, then each quasigroup
isotopic to $Q(\cdot)$ has one also. Any group of odd order has a
complete mapping, but, for example, groups of order $4k+2$ do not
contain such mappings. More interesting facts on the Brualdi
conjecture one can find in \cite{Bala, VD1, DK, Fu} and
\cite{Wan}.

If this conjecture is true, then from our results it follows that
{\it each finite quasigroup has a prolongation}.

\footnotesize{\hfill Received \ September 30, 2008

\smallskip\noindent
I.I.Deriyenko:\\ Kremenchuk State Polytechnical University,
Pervomayskaya 20, 39600 Kremenchuk, Ukraine\\
 E-mail: {ivan.deriyenko@gmail.com}\\[6pt]
W.A.Dudek:\\
 Institute of Mathematics and Computer Science,
Wroclaw University of Technology, Wybrzeze Wyspianskiego 27,
50-370 Wroclaw, Poland\\
E-mail: {dudek@im.pwr.wroc.pl} }

\end{document}